\documentclass[a4paper,11pt]{article}
\usepackage{graphicx}
\usepackage{amsmath}
\usepackage{amsfonts}
\usepackage{amssymb}
\advance\textwidth 1.2in
\advance\textheight 0.55in
\advance\oddsidemargin -0.6in
\advance\topmargin -0.3in
\usepackage{graphicx}
\usepackage{amsmath}
\usepackage{amsfonts}
\usepackage{amssymb}
\usepackage{tikz}
\usetikzlibrary{arrows}
\tikzstyle{block}=[draw opacity=0.7,line width=1.4cm]
\newtheorem{algt}{Algorithm}
\begin{document}
\title{A Hybrid Factorization Algorithm for Sparse Matrix with Mixed Precision Arithmetic}
\author{Atsushi SUZUKI
\thanks{
  RIKEN Center for Computational Science\,\ %
  Minatojima-minami-machi, Chuo, Kube Hyogo 650-0047, Japan\newline
  \hspace*{1.8em}{e-mail: {\tt Atsushi.Suzuki.aj@a.riken.jp}}\newline
  \hspace*{1.8em}submitted to proceedings of ECCOMAS Congress 2022
}}
\maketitle
% \abstract{
\begin{center}
  \begin{minipage}{5.5in}
  {\bf abstrct}\\
  A new hybrid algorithm for $LDU$-factorization for large sparse matrix combining iterative solver, which can keep the same accuracy as the classical factorization, is proposed. %
  The last Schur complement will be generated by iterative solver for multiple right-hand sides using  block GCR method with the factorization in lower precision as a preconditioner, which achieves mixed precision arithmetic, and then
  the Schur complement will be factorized in higher precision. %}
  In this algorithm, essential procedure is decomposition of the matrix into a union of moderate and hard parts, which is realized by $LDU$-factorization in lower precision with
  symmetric pivoting and threshold postponing technique.
\end{minipage}
\end{center}
\section{Introduction}
There are two kinds of demand in finding a solution of linear system with large sparse matrix in numerical simulation by using mixed precision arithmetic. %
One is for solving the system with very high condition number in numerical
simulation of complex physical model and/or with large variety of physical coefficients. %
In this case, a monolithic direct factorization solver using ``quadruple'' precision could be only feasible tool. %
However  arithmetic complexity by ``double-double'' data structure, which is a faster implementation of ``quadruple'' arithmetic using modern hardware named as fused multiply-add unit, is 25 times higher than ``double''. %
Hence, it is necessary to introduce mixed precision arithmetic with ``double'' to reduce computational complexity. 
The other is to solve rather moderate problem on forthcoming CPU with more single precision units than double. %
\par
Usages of mixed precision arithmetic in numerical linear algebra are hot
research topics and a survey paper \cite{HighamMary:2022} covers the recent developments. %
The main tool to improve or recover the accuracy of the solution obtained by lower precision either direct solver or iterative solver is the iterative refinement, which generates new right-hand side from the residual to improve the solution. %
However convergence of the refinement process depends on the condition number and it is still not easy to improve the solution for matrix with high condition number.
For some singular matrix whose condition number on the image is moderate, accurate factorization is mandatory especially to perform rank-revealing.
Our aim is to construct factorization itself, which fits to usage of mixed precision arithmetic, not to improve the accuracy of the solution.
\par
In section \ref{sec:2}, classical $LDU$-factorization with symmetric pivoting with threshold postponing is viewed and a way to decomposition of the matrix into a union of moderate and hard parts is proposed. Section~\ref{sec:3} describes a novel method to generate Schur complement matrix of the hard part by solving linear system of the moderate part with multiple right-hand sides. Section~\ref{sec:4} verifies efficiency of the proposed algorithm in accuracy and computing time by numerical examples.
\section{Factorization with symmetric pivoting for large sparse matrices}\label{sec:2}
Let us assume that the matrix $K\in\mathbb{R}^{L\times L}$ is scaled so that diagonal entries take one of $-1$, $0$ and $1$, which could be realized by a scaling with only diagonal entries $[Q]_{i\,i}=1/\textstyle{\sqrt{[\overline{K}]_{i\,i}}}$ when $[\overline{K}]_{i\,i}\neq 0$, otherwise $[Q]_{i\,i}=1$. Here the original matrix $\overline{K}$ is scaled as $Q\overline{K}Q=K$. %
\subsection{factorization with symmetric pivoting and recursive generation of Schur complement}
The $LDU$-factorization of the matrix $K$ consists of recursive generation of the Schur complement, %
\begin{align}
    \Pi^{T}\begin{bmatrix}
      K_{11} & K_{12} \\ K_{21} & K_{22} 
    \end{bmatrix}\Pi
     & =
\Pi^{T}
   \begin{bmatrix}
     K_{11} & 0 \\ K_{21} & S_{22} 
   \end{bmatrix}
   \begin{bmatrix}
I_{1} & K_{11}^{-1}K_{12}\\                              
     0 & I_{2}
   \end{bmatrix}
         \Pi\,, \label{eqn:two-by-two}\\
S_{22} &=K_{22}-K_{21}K_{11}^{-1}K_{12}\,.\label{eqn:Schur-complement}
  \end{align}
  Here pivoting strategy is a symmetric one that is expressed by using the permutation $\Pi$ with $\Pi^{T}=\Pi^{-1}$ and $\Pi$ is decomposed into a union of diagonal blocks as
  $ \Pi= \text{diag}\{\Pi_{1}\,, \Pi_{2}\}$\,.
  Indices $\Lambda=\{1,\cdots, L\}$ is decomposed into a direct sum of $\overline{\Lambda}_{1}\oplus \overline{\Lambda}_{2}$, where $\Pi_{i}$ is a one to one operation on $\overline{\Lambda}_{i}$. %
  $\overline{\Lambda}_{1}$ denotes indices of already factorized part of the matrix and further factorization will be performed for $S_{22}$ on $\Lambda_{2}$. %
Each entry of $\Pi$ is selected during the factorization. %
Let us suppose that $k\times k$ sub-matrix $K_{11}$ is already factorized and $\{\Pi(1), \cdots, \Pi(k)\}$ are obtained and in the rest of matrix $(L-k)\times (L-k)$ sub-matrix $K_{22}$ needs to be factorized. %
We find the maximum absolute value in the diagonal entries of $K_{22}$, label $k_{0}$ as the index of such entry, define $\Pi(k+1)=k_{0}$, and set $d=[K]_{k_{0}\,k_{0}}$. %
The rows and columns of $K$ with indices $k+1$ % <= $k$ :: 27 Jul.2022
and $k_{0}$ are exchanged each other and the result is stored as $\widetilde K$. Schur complement matrix $S_{22}'$, whose size is $(L-k-1)\times (L-k-1)$, is calculated by rank-$1$ update with weight $1/d$, $S_{22}'={\widetilde K}_{22}'-[{\widetilde K}_{22}]_{\downarrow\,1}d^{-1}[{\widetilde K}_{22}]_{1\,\rightarrow}$. This is the essential operation of $LDU$-factorization. %
\par
In practical computation of large sparse matrix,  $\Lambda_{1}$ will be created by following the elimination tree of the nested-dissection ordering with threshold postponing with user defined parameter $\tau$. %
  \subsection{nested-dissection ordering and threshold postponing}
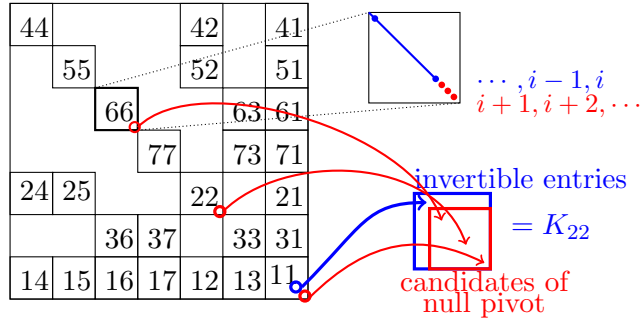
\begin{figure}[tb]
  \centering
  \begin{tikzpicture}[scale=0.8]
  \draw (6.5,3.25) rectangle (8.,4.75);
  \draw[blue,thick] (6.5,4.75) -- (7.6,3.65);
  \fill[blue] (6.6,4.65) circle (0.33ex);
  \fill[blue] (7.6,3.65) circle (0.33ex);
  \fill[red] (7.7,3.55) circle (0.33ex);    
  \fill[red] (7.8,3.45) circle (0.33ex);  
  \fill[red] (7.9,3.35) circle (0.33ex);
  \node[blue,right] at (8.1,3.6) {\small $\cdots, i-1, i$};
  \node[red,right] at (8.1,3.2) {\small $i+1, i+2, \cdots $};
  \draw (0.6,0) rectangle (5.5,4.9);
  \draw (0.6,4.2) rectangle (1.3,4.9);
  \node at (1.0,4.5) {44};
  \draw (1.3,3.5) rectangle (2.0,4.2);
  \node at (1.7,3.8) {55};
  \draw (2.0,2.8) rectangle (2.7,3.5);
  \node at (2.4,3.1) {66};
  \draw (2.7,2.1) rectangle (3.4,2.8);
  \node at (3.1,2.4) {77};
  \draw (3.4,1.4) rectangle (4.1,2.1);
  \node at (3.8,1.7) {22};
  \draw (4.1,0.7) rectangle (4.8,1.4);
  \node at (4.5,1.0) {33};
  \draw (4.8,0.0) rectangle (5.5,0.7);
  \node at (5.1,0.4) {11};
  \draw (4.8,0.7) rectangle (5.5,1.4);
  \node at (5.2,1.0) {31};
  \draw (4.8,1.4) rectangle (5.5,2.1);
  \node at (5.2,1.7) {21};
  \draw (4.8,2.1) rectangle (5.5,2.8);
  \node at (5.2,2.4) {71};
  \draw (4.8,2.8) rectangle (5.5,3.5);
  \node at (5.2,3.1) {61};
  \draw (4.8,3.5) rectangle (5.5,4.2);
  \node at (5.2,3.8) {51};
   \draw (4.8,4.2) rectangle (5.5,4.9);
  \node at (5.2,4.5) {41};
  \draw (4.1,2.1) rectangle (4.8,2.8);
  \node at (4.5,2.4) {73};
  \draw (4.1,2.8) rectangle (4.8,3.5);
  \node at (4.5,3.1) {63};
  \draw (3.4,3.5) rectangle (4.1,4.2);
  \node at (3.8,3.8) {52};
  \draw (3.4,4.2) rectangle (4.1,4.9);
  \node at (3.8,4.5) {42};
  \draw (0.6,0.0) rectangle (1.3,0.7);  
  \node at (1.0,0.3) {14};
  \draw (1.3,0.0) rectangle (2.0,0.7);    
  \node at (1.7,0.3) {15};
  \draw (2.0,0.0) rectangle (2.7,0.7);    
  \node at (2.4,0.3) {16};
  \draw (2.7,0.0) rectangle (3.4,0.7);    
  \node at (3.1,0.3) {17};
  \draw (3.4,0.0) rectangle (4.1,0.7);    
  \node at (3.8,0.3) {12};
  \draw (4.1,0.0) rectangle (4.8,0.7);  
  \node at (4.5,0.3) {13};
  \draw (2.0,0.7) rectangle (2.7,1.4);    
  \node at (2.4,1.0) {36};
  \draw (2.7,0.7) rectangle (3.4,1.4);      
  \node at (3.1,1.0) {37};
  \draw (0.6,1.4) rectangle (1.3,2.1);
  \node at (1.0, 1.8) {24};
  \draw (1.3,1.4) rectangle (2.0,2.1);
  \node at (1.7,1.8) {25};
  
  \draw[thick] (2.0,2.8) rectangle (2.7,3.5);
  \draw[densely dotted] (2.0,3.5) -- (6.5,4.75);
  \draw[densely dotted] (2.7,2.8) -- (8,3.25);
  \draw[blue,very thick](7.25,0.5) rectangle (8.5,1.75);
  \node[blue] at (9.5,1.2) {$=K_{22}$};
  \node[blue] at (9,2) {invertible entries};
  \node[red] at (8.4,0.3) {candidates of};
   \node[red] at (8.4,-0.1) {null pivot};
   \draw[red,very thick] (7.5,0.5) rectangle (8.5,1.5);
   \draw [->,red, thick] (5.52,0.02) .. controls +(1,1) and + (-1,0.5) .. (8.4,0.6);
   \draw[red, very thick] (5.45,0.05) circle (0.5ex);
   \draw[blue, very thick] (5.3,0.2) circle (0.5ex);
   \draw [->,blue, very thick] (5.4,0.2) .. controls +(1,1) and + (-1,0) .. (7.45,1.6);
   \draw[red, thick] (2.65,2.85) circle (0.5ex);
   \draw [->,red, thick] (2.7,2.85) .. controls +(1,1) and + (-0.5,1.5) .. (7.7,1.3);
   \draw[red, very thick] (4.05,1.45) circle (0.5ex);
   \draw [->,red, thick] (4.1,1.45) .. controls +(1,1) and + (-0.5,1.5) .. (8.1,0.9);
 \end{tikzpicture}
  \caption{symmetric pivoting with threshold postponing following nested-dissection ordering}
  \label{fig:threshold-postponing}
\end{figure}

  Let us introduce the nested-dissection ordering~\cite{Geroge:1977} and suppose that $\Lambda$ is decomposed into $2^{m}-1$ sub-indices with $m$-level bi-section tree, $\Lambda=\bigoplus_{1\leq \ell \leq 2^{m}-1}\Lambda_{\ell}$, where $k$-th level contains $2^{k-1}$ sub-index sets $\{\Lambda_{\ell}\}$. %
  For stability of the factorization, a given threshold parameter $\tau > 0$ is introduced to perform postponing of factorization. %
  During the $LDU$-factorization of the sub-matrix with index $\Lambda_{j}$, if the ratio in absolute value of successive diagonal entries becomes smaller than $\tau$, i.e., $|[K]_{i+1\,i+1}/[K]_{i\,i}|< \tau$, then the lower block of the matrix is not factorized. %
  The index will be decomposed as  $\Lambda_{j}=\{j_{1},\cdots, j_{i}\}\oplus\{j_{i+1},\cdots, j_{M}\}=\widetilde{\Lambda_{j}}\oplus\widehat{\Lambda_{j}}$,
  where $\widetilde{\Lambda_{j}}$ is set of the indices for the factorized part.
  For $m$-level bisection tree of the nested-dissection ordering,
 an $LDU$-factorizable part with $1\times 1$ pivot is collected as
$ \bigoplus_{1\leq \ell\leq 2^{m}-1}\widetilde{\Lambda_{\ell}}$
and $\Lambda_{0}=\bigoplus_{1\leq \ell\leq 2^{m}-1}\widehat{\Lambda_{\ell}}$ for the postponed entries. %
At the end of all threshold factorization following the elimination tree,
we will again apply the $LDU$-factorization to the last Schur complement with indices $\Lambda_{0}$ and will obtain $\Lambda_{0}=\widetilde{\Lambda_{0}}\oplus \widehat{\Lambda_{0}}$. %
Here we will enlarge $\widehat{\Lambda_{0}}$ with $n$ entries by moving the last entries of $\widetilde{\Lambda_{0}}$ to ensure $S_{22}$ has an image space, which contributes to comparison between zero eigenvalues and nonzero ones. %
Usually we take $n=4$ entries~\cite{SuzukiRoux:2014}. %
By this process, the very last Schur complement matrix $K_{22}$ in \eqref{eqn:two-by-two} collecting postponed pivots $\Lambda_{2}:=\bigoplus_{0\leq\ell\leq 2^{m}-1}\widehat{\Lambda_{\ell}}$ will have large condition number or singular for the case that the original matrix is not invertible, and may contain $2\times 2$ pivoting entries when $K_{11}$ is not definite. %
On the contrary, $K_{11}$ in \eqref{eqn:two-by-two} with index $\Lambda_{1}:=\bigoplus_{0\leq\ell\leq 2^{m}-1}\widehat{\Lambda_{\ell}}$ has moderate condition number that can be factorized with appropriate permutation.
Figure~\ref{fig:threshold-postponing} shows schematic explanation of the threshold postponing with $3$-level bisection tree in nested-dissection ordering.
  \section{Hybrid factorization algorithm}\label{sec:3}
  A new algorithm is constructed by replacing the solution of the linear system with multiple right-hand sides (RHSs), $K_{11}X_{12}=K_{12}$ of the first block in \eqref{eqn:two-by-two} by an iterative solver. %
  Factorization for $K\in \mathbb{R}^{L\times L}$ with user-defined threshold $\tau$ is performed in a hybrid way as follows.
  \begin{algt} \label{algt:1} $LDU$-factorization with internal iterative solver \\
    \rm
1.\ \ factorize $K$ with threshold postponing in lower precision and extract moderate part $K_{11}$ with 
    finding indices $\Lambda_{1}$ and permutation $\Pi_{1}$
   with $N=\#\Lambda_{1}$\,.\\
2.\ \ decompose matrix $K$ into $2\times 2$ blocks as $\begin{bmatrix}
     K_{11} & K_{12} \\ K_{21} & K_{22} 
   \end{bmatrix}$ 
   with $K_{11}\in\mathbb{R}^{N\times N}$, $K_{12}\in\mathbb{R}^{N\times M}$, $K_{21}\in\mathbb{R}^{M\times N}$, and $K_{22}\in\mathbb{R}^{M\times M}$\,.\\
   3.\ \ find solution $X_{12}$ satisfying $K_{11}X_{12}=K_{12}$ by an iterative solver using the $LDU$-factoriza\-tion of $K_{11}$ with permutation $\Pi_{1}$ in lower precision as preconditioner\,.\\
4.\ \ construct the Schur complement $S_{22}:=K_{22}-K_{21}X_{12}$ in higher precision\,.\\
5.\ \ factorize $S_{22}$ in higher precision with finding a symmetric pivoting expressed by the permutation $\Pi_{2}$ that may contain $2\times 2$ entries\,.
\end{algt}
We can utilize this solution in lower precision as a preconditioner for the iterative solver in higher precision.
\par
If the condition number of $K_{11}$ is in the range of the maximum floating digits of the lower precision, the solution $X_{12}$ of $K_{11}X_{12}=K_{12}$ will be obtained very accurately with the residual closed to the machine epsilon of the higher precision.
Therefore, even for the case that the original matrix is singular, the Schur complement $K_{22}$ is well constructed without large perturbation during the operation of $K_{11}^{-1}$ and kernel detection for rank-reveling~\cite{SuzukiRoux:2014} works well as the original full factorization algorithm. %
Solution procedure with forward and backward substitutions are performed block-wisely for the linear system after applied permutations $\Pi = \text{diag}[\Pi_{1}, \Pi_{2}]$. The following algorithm describes a procedure to find $[x_{1}^{T}, x_{2}^{T}]^{T}$ satisfying
\begin{equation}
  \begin{bmatrix}
    K_{11} & K_{12} \\
    K_{21} & K_{22}
  \end{bmatrix}
  \begin{bmatrix}
    x_{1} \\ x_{2}
  \end{bmatrix}
  =
    \begin{bmatrix}
    b_{1} \\ b_{2}
  \end{bmatrix}\,.\label{eqn:linear-system-in-two-blocks}
\end{equation}
\begin{algt} \label{algt:2} forward/backward substitutions for hybrid factorization\\
  \rm
1. \ solve $K_{11}\,y_{1}=b_{1}$ by the same iterative solver in Algorithm~\ref{algt:1} in higher precision\,.\\
2. \ compute $y_{2}=b_{2}-K_{21}y_{1}$\,.\\
3. \ solve $S_{22}\,x_{2}=y_{2}$ by forward and backward substitutions in higher precision\,.\\
4. \ update $x_{1}=y_{1}-X_{12}x_{2}$ with $X_{12}$ is computed in Algorithm~\ref{algt:1}\,.
\end{algt}
Solution of the moderate part of the liner system $K_{11}y_{1}=b_{1}$ is obtained by the iterative solver within given accuracy to specify the convergence. %
Thanks to preconditioner by the $LDU$-factorization in lower precision for $K_{11}$ with moderate condition number, we can expect enough accuracy of the solution by the preconditioned iterative solver closed to one by the direct solver. %
\subsection{iterative solver for solution of linear system with multiple RHSs}
There are two kinds of solver for the system $K_{11}X_{12}=K_{12}$.
For simplicity, we first describe algorithms for linear system with single RHS, $A\,x = b$ where $A\in\mathbb{R}^{N\times N}$ is invertible and it stands for $K_{11}$ and $b\in\mathbb{R}^{N}$ will be one of the column vector of $K_{12}$.
Let us denote $\widehat{Q}\widehat{x}= \widehat{b}$ as the linear system in lower precision with solution $\widehat{x}$ for the RHS, $\widehat{b}$, which is converted from the given data $b$ in higher precision by the floating point casting operation.
\subsection{iterative refinement}
The iterative refinement is a classical method to improve the accuracy of the linear system. %
For mixed precision arithmetic, solution of the linear system is found in lower precision but calculation of the residual is performed in higher precision. %
Therefore we can expect iterative renfinement will converge with higher accuracy.
\begin{algt}\label{algt:3} iterative refinement to improve solution in lower precision\\
  \rm
1. \ find $\widehat{e_{0}}$ satisfying $\widehat Q\widehat{e_{0}}=\widehat{b}$ in lower precision\,.\\
2. \ convert  $x_{0}\leftarrow\widehat{e_{0}}$ from lower precision to higher precision\,.\\
3. \ compute residual $r_{0}=b-A\,x_{0}$ of solution\,.\\
4. \ loop $n=0, 1, 2, \cdots$\\
\mbox{}\quad 4a. \ truncate $\widehat{r_{n}}\leftarrow r_{n}$, higher precision data to lower preicision\,.\\
\mbox{}\quad 4b. \ find $\widehat{e_{n}}$ satisfying $\widehat Q\widehat{e_{n}}=\widehat{r_{n}}$ in lower precision\,.\\
\mbox{}\quad 4c. \ update solution $x_{n+1}=x_{n}+\widehat{e_{n}}$ by adding lower precision data\,.\\
\mbox{}\quad 4d. \ compute residual $r_{n+1}=b-A\,x_{n+1}$ of solution\,.
\end{algt}
Here in the procedure for updating $(n+1)$-th solution, addition of lower precision data can be performed directly without preparation of a working vector $e_{n}$ in higher precision converting form lower precision $\widehat{e_{n}}$. %
To make clear of the role of the preconditioner $Q$ in higher precision,
let $Q\,e_{n}=r_{n}$ denote the operation to find $\widehat{e_{n}}$ satisfying $\widehat{Q}\widehat{e_{n}}=\widehat{r_{n}}$ in lower preicision for given data  $r_{n}$ that is converted to $\widehat{r_{n}}$ and up-converting $\widehat{e_{n}}$ in lower to $e_{n}$. %
Calculation of the residual $r_{1}$ in the first step is viewed as following using the assumption that $A$ is invertible, %
  \begin{align*}
    r_{1}&=b-A\,x_{1}=b-A(x_{0}+e_{1})=r_{0}-Ae_{1}=r_{0}-AQ^{-1}r_{0}=(I-AQ^{-1})r_{0}\,,\\
    x_{1}&=x_{0}+A^{-1}(I-(I-AQ^{-1}))r_{0}=x_{0}+A^{-1}AQ^{-1}r_{0}=x_{0}+Q^{-1}r_{0}\,. %    
  \end{align*}
  By the same argument, residual and solution at $n$-th step are obtained as
  \begin{align*}
    r_{n}&=(I-AQ^{-1})^{n}r_{0}\\
    x_{n}&=x_{0}+A^{-1}(I-(I-A\,Q^{-1})^{n})r_{0}\\
         &=x_{0}+\dbinom{n-1}{0}Q^{-1}r_{0}
           -\dbinom{n-1}{1}(Q^{-1}A)Q^{-1}r_{0}+ \cdots
           +(-1)^{(n-1)}\dbinom{n-1}{n-1}(Q^{-1}A)Q^{n-1}r_{0}
  \end{align*}
  using the binomial expansion. %
  We conclude that iterative refinement will find solution $x_{0}$ in a Krylov subspace with $Q^{-1}A$ and $Q^{-1}r_{0}$
  \begin{equation*}
    x_{n}\in x_{0}+\text{span}[Q^{-1}r_{0},
    (Q^{-1}A)Q^{-1}r_{0},
    (Q^{-1}A)^{2}Q^{-1}r_{0}, \cdots,
    (Q^{-1}A)^{n-1}Q^{-1}r_{0}]\,. %
  \end{equation*}
\subsection{preconditioned GCR method for single RHS}
The iterative refinement procedure to improve accuracy of the solution obtained by lower precision arithmtic can be viewed as an iterative process to find solution in the Krylov subspace of preconditioned matrix with fixed coefficient for linear combination. %
The coefficient by the binomial expansion is not optimal and we can use the standard proceudre of Krylov subspace solver family. %
The most easiest method in implementation is Generalized Conjugate Residual (GCR) method~\cite{Saad:2003} and it is also closed to the iterative refinement procedure with further approximation. %
A preconditioned GCR method for the linear system $A\,x=b$ with $A\in \mathbb{R}^{N\times N}$ and $b\in \mathbb{R}^{N}$ by using solution $Q\,y=f$ with $Q\in \mathbb{R}^{N\times N}$ and $f\in \mathbb{R}^{N}$ in lower precision as a right preconditioner
\begin{equation*}
  AQ^{-1}Q\,x=AQ^{-1}\widetilde{x}=b
\end{equation*}
is given as Algorithm~\ref{algt:4}.
When $A$ has moderate condition number, the solution process by $LDU$-factorization of $A$ in lower precision that is expressed as $Q^{-1}$ is well performed and $AQ^{-1}$ is very closed the identity matrix, $AQ^{-1}\simeq I_{N}$. %
In practice the following right preconditioned GCR converges in few iterations, which is rather natural consequence by selection of moderate part of the matrix $K_{11}$ using threshold postponing in lower precision. %
An example of convergence history will be shown in Section~\ref{sec:GCR-conv}.
\noindent
\begin{algt}\label{algt:4}preconditioned GCR method\\
  \rm
find $x_{0}$ satisfying $Q\,x_{0}=b$\\
$r_{0}=b-A\,x_{0}$\\
$w_{0}=Q^{-1}r_{0}$\\
$p_{0}=w_{0}$\\
loop $n = 0, 1, 2, \cdots$\\
\mbox{}\quad $\alpha_{n}:=\dfrac{(r_{n}, Ap_{n})}{(A\,p_{n}, A\,p_{n})}=\dfrac{(r_{n}, AQ^{-1}\widetilde{p_{n}})}{(AQ^{-1}\widetilde{p_{n}}, AQ^{-1}\widetilde{p_{n}})}$\\
\mbox{}\quad $x_{n+1}:=x_{n}+\alpha_{n}p_{n}$\\
\mbox{}\quad $r_{n+1}:=r_{n}-\alpha_{n}A\,p_{n}$\\
\mbox{}\quad find $w_{n+1}$ satisfying $Q\,w_{n+1}=r_{n+1}$\\
\mbox{}\quad for $0\leq m\leq n$\\
\mbox{}\quad\quad $\beta_{m\,n}:=-\dfrac{(A\,w_{n+1}, A\,p_{n})}{(A\,p_{m}, A\,p_{m})}=
-\dfrac{(AQ^{-1}r_{n+1}, AQ^{-1}\widetilde{p_{n}})}{(AQ^{-1}\widetilde{p_{m}}, AQ^{-1}\widetilde{p_{m}})}$\\
\mbox{}\quad $p_{n+1}:=w_{n+1}+\displaystyle{\sum_{m=0}^{n}}\beta_{m\,n}p_{m}$
\end{algt}
In practical computation, to avoid two times multiplication of $A$ to $p_{n}$ and to $w_{n+1}$, $A\,p_{n}$ is stored as $q_{n}$ and is updated in the same manner as $p_{n+1}$ using a new vector $z_{n+1}:=A\,w_{n+1}$, which results in
$q_{n+1}:=z_{n+1}+\sum_{m=0}^{n}\beta_{m\,n}q_{m}$. %
% This treatment will be directly described in Algorithm~\ref{algt:5} for multiple RHSs. %
Here $A\,w_{n+1}$ is performed by the SpMV (Sparse Matrix-Vector multiplication) opertaion. %
\par
We can see the residual at $n$-step belongs to $(n+1)$-dimensional Krylov subspace, $r_{n}\in K_{n+1}(r_{0}, AQ^{-1})$ and approximate solution $\widetilde{x}_{n}$  of the right preconditioned system $AQ^{-1}\widetilde x=b$ is found as
$\widetilde{x}_{n}\in \widetilde{x}_{0}+K_{n}(r_{0}, AQ^{-1})$ with the initial approximation $Q\,x_{0}=\widetilde{x}_{0}$. By defintion of $\widetilde{x}_{n}=Q\,x_{n}$, approximation of the linear system is written as
${x}_{n}\in {x}_{0}+K_{n}(Q^{-1}r_{0}, Q^{-1}A)$. %
\par
Orthogonality on residual, for all $y\in K_{n}(r_{0},AQ^{-1})$, it holds
that $(r_{n+1}, AQ^{-1}y)=0$ and search vectors 
$(AQ^{-1}\widetilde{p_{m}}, AQ^{-1}\widetilde{p_{n}})=0$ for $m\neq n$. Both are verified by induction.
Since we can assume $AQ^{-1}\simeq I_{N}$, we get
\begin{equation*}
  \alpha_{n}=\dfrac{(r_{n}, AQ^{-1}\widetilde{p_{n}})}{(AQ^{-1}\widetilde{p_{n}}, AQ^{-1}\widetilde{p_{n}})}\simeq
  \dfrac{(r_{n}, \widetilde{p_{n}})}{(\widetilde{p_{n}}, \widetilde{p_{n}})}\quad \text{ and }
  \beta_{m\,n}\simeq - \dfrac{(r_{n+1}, \widetilde{p_{n}})}{(\widetilde{p_{m}},\widetilde{p_{m}})}\,. %
\end{equation*}
If we could approximate $\alpha_{n}\simeq 1$ and $\beta_{m\,n}\simeq 0$ with $0\leq m \leq n$ up to $n$-th step, we will have
\begin{equation*}
  x_{n+1}\simeq x_{n}+p_{n},\quad r_{n+1}\simeq r_{n}-A\,p_{n},\quad \text{ and }\quad p_{n}\simeq Q^{-1}r_{n+1},
\end{equation*}
which leads to the same procedure of the iteartive refinement. %
\par
We can expect that the preconditioned GCR method converges faster than the the iterative refinement thanks to better combination of coefficients for Krylov subspace basis to achieve the Galerkin orthogonality.
\subsection{preconditioned block GCR method for multiple RHSs}
In practice, the linear system consists of multiple RHSs, $K_{12}$ and then it will be more efficient to use block GCR method for multiple RHS, because dimension of Krylov subspace in block version is much larger than one for single RHS.
Comparison of convergence of single and multiple RHSs versions will be illustrated in Section~\ref{sec:GCR-conv}. %
\par
A preconditioned block GCR method for the linear system for multiple RHSs, $A\,[x^{(1)}, \cdots,$\\$ x^{(M)}]=[b^{(1)}, \cdots, b^{(M)}]$ by using solution $Q\,y=f$ in lower precision as a right preconditioner is obtained by introducing multiplication of matrix $A$ to $[w_{n+1}^{(1)},\cdots, w_{n+1}^{(M)}]$, which is called SpMM (Sparse Matrix-Matrix multiplication) operation. %
Determination of magnitude of the search vector $\alpha_{n}$ in updating procedure is replaced by operation by $M\times M$ matrices as
\begin{align*}
[{\cal M}_{n}]_{\ell, k}&:=(q_{n}^{(k)}, q_{n}^{(\ell)})\ \text{ and  }\ %
[{\cal A}_{n}]_{\ell, k}:=(r_{n}^{(k)}, q_{n}^{(\ell)})\ \text{ for } \ 1\leq \ell, k\leq M\\
[x_{n+1}^{(1)},\cdots, x_{n+1}^{(M)}]&:=[x_{n}^{(1)}, \cdots, x_{n}^{(M)}]+[{p_{n}}^{(1)}, \cdots, {p_{n}}^{(M)}]{\cal M}_{n}^{-1}{\cal A}_{n}\,.
\end{align*}
During the iteration, $[q_{n+1}^{(1)},\cdots, q_{n+1}^{(M)}]$ are updating with keeping identity to multiplication of $A$ to $[p_{n+1}^{(1)},\cdots, p_{n+1}^{(M)}]$. %
Both sets of vectors are calculated from linear combination expressed by ${\cal M}_{n}^{-1}{\cal B}_{m\,n}$ with
$[{\cal B}_{m\,n}]_{\ell, k}:=-(z_{n+1}^{(k)}, q_{n}^{(\ell)})$. %
We call this iteration as {\bf Algorithm 5}.
When all matrices ${\cal M}_{n}$ at $n$-th iteration are invertible, the approximate solution $x_{n}^{(k)}$ with $1\leq k\leq M$ is found as
$x_{n}^{(k)}\in x_{0}^{(k)}+K_{n}(Q^{-1}[r_{0}^{(1)},\cdots, r_{0}^{(M)}], Q^{-1}A)$
and 
$(r_{n}^{(k)}, AQ^{-1}y)=0$ for any vector $y \in K_{n}([r_{0}^{(1)},\cdots, r_{0}^{(M)}], AQ^{-1})$.
\subsection{convergence history of preconditioned GCR method}\label{sec:GCR-conv}
\begin{figure}[tb]
  \centering
  \includegraphics[width=3.5in]{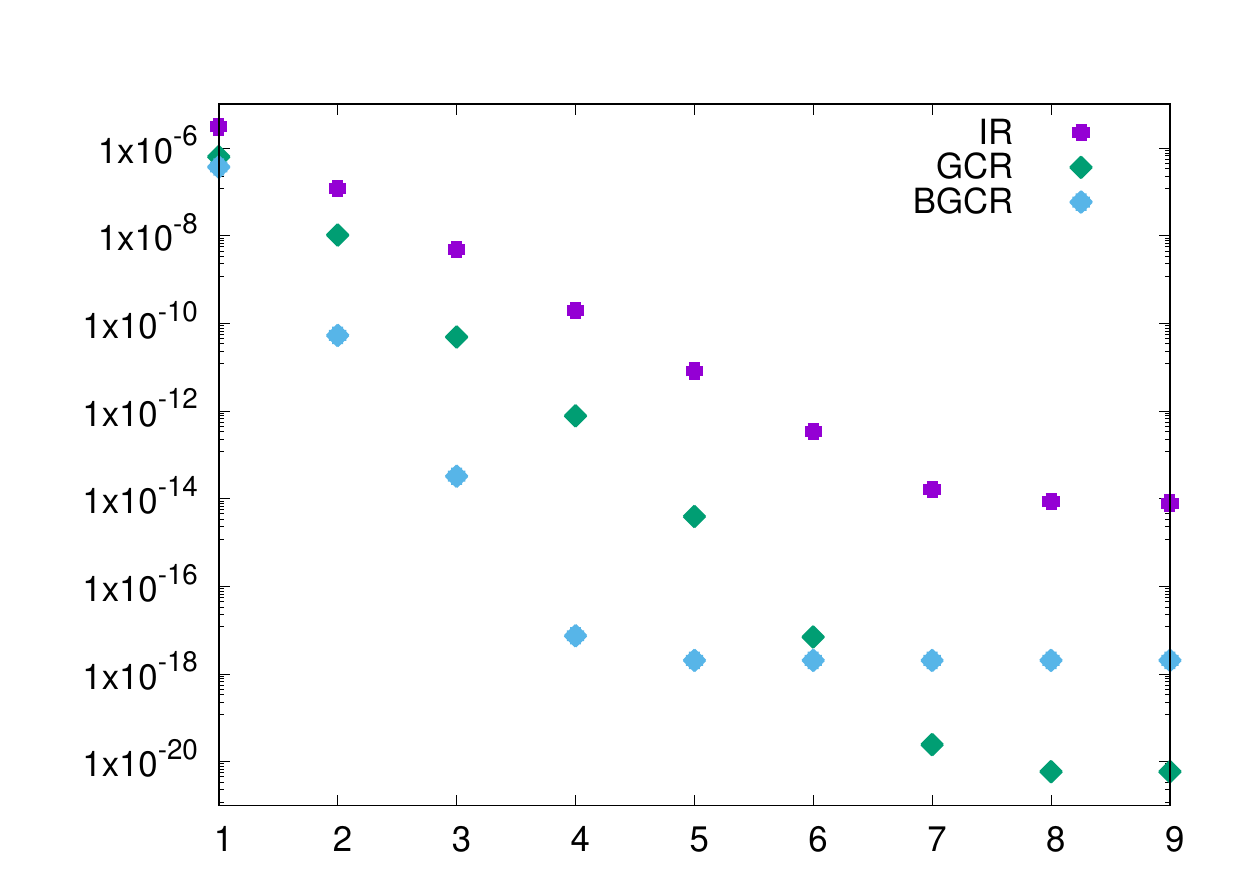}
  \caption{convergence history of IR/GCR/BGCR}
  \label{fig:conveGCR}
\end{figure}

Figure~\ref{fig:conveGCR} shows convergence history to find a solution of $X_{12}$ satisfying $K_{11}X_{12}=K_{12}$ with matrix in Section~\ref{sec:numerical-Stokes}. %
Sixteen entries are postponed during factorization of $K$ by single precision arithmetic and $K_{12}$ consists of 16 column vectors with $N=374,520$. %
Each iteration drawn by colors, e.g., purple for iterative refinement (IR): Algorithm~\ref{algt:3}, green for preconditioned GCR: Algorithm~\ref{algt:4} , and light blue for preconditioned block GCR (BGCR): Algorithm~5, shows convergence of the first column vector of sixteen RHSs. %
We can see convergence of IR is slower than other GCR solvers and block GCR converges after 5 iterations to the machine epsilon of the double precision. %
\section{Numerical examples}\label{sec:4}
In this section, feasibility of hybrid factorization and solution Algorithms~\ref{algt:1} and \ref{algt:2} with inner iterative solver by Algorithm 5 using preconditioner in lower precision solution, will be demonstrated by three sparse matrices. %
\subsection{symmetric matrix from matrix market}\label{sec:matrix-market}
The matrix {\tt phb1HYS} from the matrix market database is symmetric and has $36,414$ total number of unknowns and $4,344,765$ originated from protein data bank and was used in evaluation of SpMV performance in \cite{WilliamsOlikerVuducShalfYelickDemmel}.
By setting pivot threshold $\tau=0.05$, 17 entries are excluded by factorization in single precision. %
By this decomposition, Algorithm~\ref{algt:1} is performed with $N=\#\Lambda_{1}=36,397$ and $M=\#\Lambda_{2}=17$. %
The maximum and minimum eigenvalues are calculated by the power method as
$\lambda_{\max}(K_{11})=1.0002$, $\lambda_{\min}(K)=3.8109\times 10^{-4}$ and then
$\kappa(K_{11})=2.6448\times 10^{3}$.
Since the eigenvalue distribution of the whole matrix has significant jumps in the smallest, $\lambda_{\min}(K)=9.9876\times 10^{-10}$,
$\lambda_{\max}(K)=34.854$ and then $\kappa(K)=3.4862\times 10^{11}$.
The single precision arithmetic can factorize the moderate part $K_{11}$ but accuracy of the solution of the whole system with $K$ is very poor. %
Since the matrix {\tt phb1HYS} is symmetric, $LDL^{T}$-factorization is performed instead of $LDU$-factorization, but block GCR : Algorithm~5 is used for the system with $K_{11}$ and 17 RHSs.
\subsection{finite element matrix from incompressible flow problem}\label{sec:numerical-Stokes}
\begin{figure}[tb]
  \centering
  \includegraphics[width=70mm]{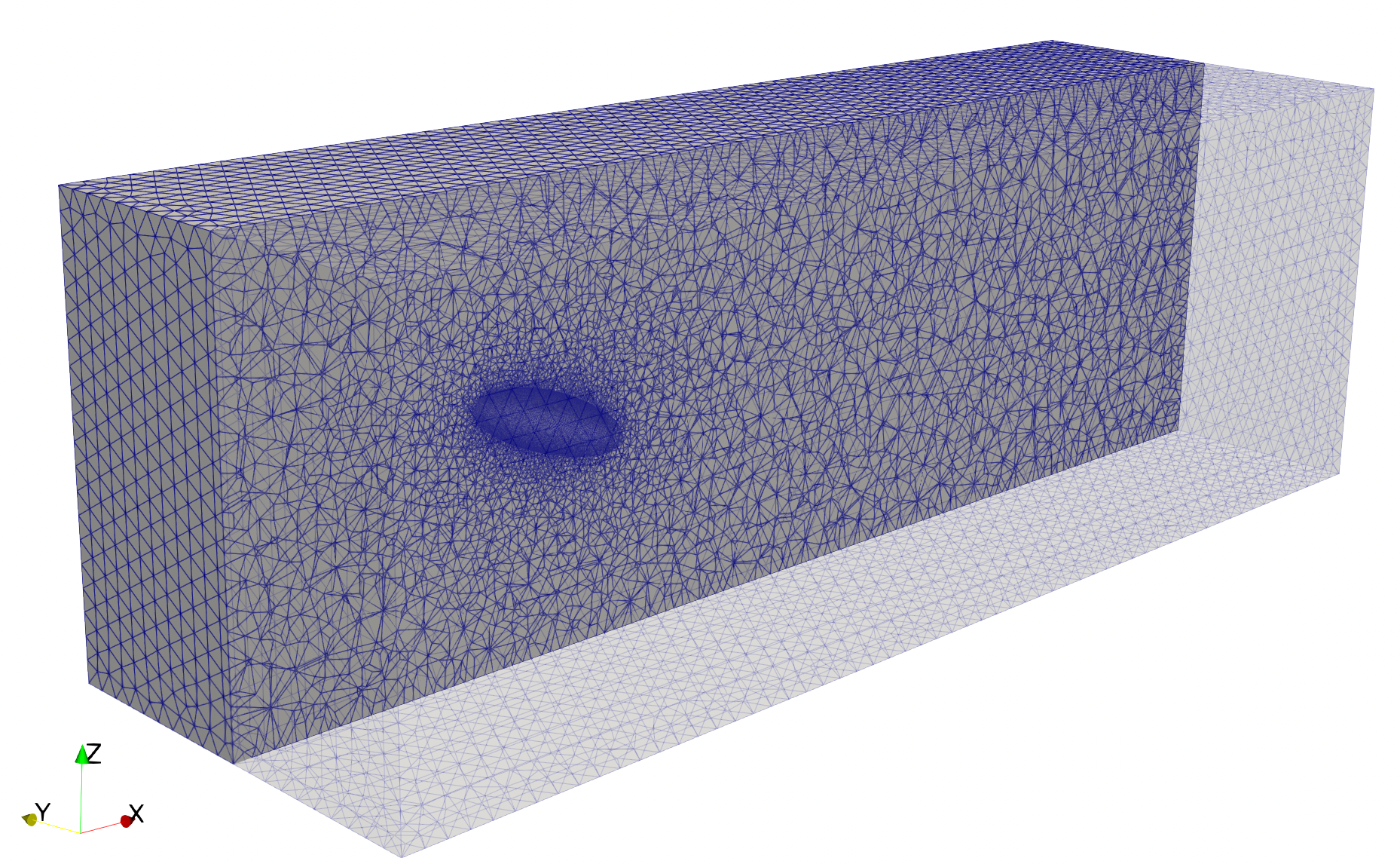}
  \caption{finite element mesh decomposition of the box domain excluding an ellipsoid}
  \label{fig:flow-excluding-ellipsoid}
\end{figure}

Let $\Omega$ be the flow region consisting of a box domain excluding an ellipsoid. %
The size of the box is $6\times 2\times 2$ and diameters of the ellipsoid is $(0.75, 0.5, 0.25)$ with 20 degree slanted. %
Finite element mesh decomposition $\overline{\Omega}=\bigcup_{e}e$ with tetrahedra $\{e\}$, whose diameter is denoted as $h_{e}=|e|$, is depicted in Figure~\ref{fig:flow-excluding-ellipsoid}. %
Here nonuniform mesh subdivision is used and smallest mesh size $h=0.01$ on the surface of the ellipsoid and largest mesh size $h=0.1$ on the inlet boundary.
\par
We consider a Stokes problem to find the velocity $u$ and the pressure $p$, $-\nabla\cdot D(u)+\nabla p=0$ and $\nabla\cdot u=0$ with full homogeneous Neumann data on the all boundary surfaces. Here $D(u)$ denotes the strain rate tensor $D(u)=(\nabla u+(\nabla u)^{T})/2$.
This boundary condition $2D(u)n-n\,p=0$ with outer normal $n$ is not physical one, but a domain decomposition method with an artificial boundary condition like FETI method~\cite{FarhatRoux:1994} leads to a floating sub-problem with such kind of full Neumann boundaries. %
Finite element matrix is obtained by discretization of the weak formulation with P1/P1 elements and a stabilization parameter $\delta$ that is set as $0.01$ in this example,
\begin{equation}
\overline{K} =
  \begin{bmatrix}
    A & B^{T} \\-B & \delta D
  \end{bmatrix}\,.\label{eqn:stiffness-Stokes}
\end{equation}
Each block is defined using finite element basis functions $\{\varphi_{i}\}$ for the velocity unknown and $\{\psi_{i}\}$ for the pressure unknown,
\begin{equation*}
  [A]_{i\,j}=\int_{\Omega} D(\varphi_{j}):D(\varphi_{i}),\quad %
  [B]_{i\,j}=-\int_{\Omega}\nabla\cdot\varphi_{j}\psi_{i},\quad %
  [D]_{i\,j}=\sum_{e}\int_{e}h_{e}^{2}\nabla\psi_{j}\cdot\nabla\psi_{i}\,.
\end{equation*}
The stiffness matrix $\overline{K}$ with total numbers of unknowns $374,536$ and nonzeros $21,146,848$ is generated from tetrahedral mesh decomposition in Figure~\ref{fig:flow-excluding-ellipsoid} by FreeFEM software package~\cite{Hecht:2002}. %
We call this matrix as {\tt stokes}. %
For constant vectors $a$ and $b$ whose total degrees of freedom is 6, we have $D(a\times x+b)=0$ and $\nabla\cdot(a\times x+b)=0$, where $a\times x+b$ called rigid body modes. Then it is clear that the stiffness matrix $\overline{K}$ defined in \eqref{eqn:stiffness-Stokes} has six dimensional kernel.
%% mistake for \lambda_\min{\text{Im}(K)} and \lambda_{K_{11}} :: 27 Jul.2022
The maximum and minimum eigenvalues are calculated by the power method as
$\lambda_{\max}(K)=2.5003$, $\lambda_{\min}(K|_{\text{Im}(A)})=8.0038\times 10^{-8}$,  and $\lambda_{\min}(K)=1.1437\times 10^{-19}$. %
Here the operator $K|_{\text{Im}(K)}$ is one-to-one in $\text{Im}(K)$, i.e., the orthogonal complement of $\text{Ker}(K)$, which is obtained numerically by $LDU$-factorization procedure.
The condition number of $K$ is $\kappa(K|_{\text{Im}(K)})=3.1240\times 10^{7}$ when it is restricted on $\text{Im}(K)$ and $\kappa(K)=2.1863\times 10^{19}$ for all unknowns.
\par
Applying single precision arithmetic with pivot threshold $\tau=0.75$, 16 entries are postponed and the Schur complement $S_{22}$, whose size is $16\times 16$, has six dimensional kernel. %
By this decomposition, $\lambda_{\max}(K_{11})=1.0000$, $\lambda_{\min}(K_{11})=1.6677\times 10^{-8}$,  and $\kappa(K_{11})=5.9962\times 10^{7}$.
\subsection{finite element matrix from semi-conductor  problem}\label{sec:semi-conductor}
The semi-conductor problem is mathematically modeled by the drift-diffusion equations with the electrostatic potential $\varphi$, the electron density $n$, and the hole density $p$. By introducing Slotboom variables for the electron density $\eta=e^{-\varphi}n$ and for the hole density $\xi=e^{\varphi}p$, the drift and diffusion terms are combined into a single term and the following nonlinear system~\cite{Brezzietal:2005} is obtained as
\begin{equation}
  -\nabla\cdot(\lambda^{2}\nabla\varphi)=e^{-\varphi}\xi-
  e^{\varphi}\eta+C(x),\quad -\nabla (e^{\varphi}\nabla\eta)=0,\quad
\nabla (e^{-\varphi}\nabla\xi)=0\label{eqn:drift-diffusion-nonlinear},
\end{equation}
where $\lambda$ denotes Debye length, and $C(x)$ is a given function to represent doping density of $N$-rich or $P$-rich material.
Here we consider two dimensional problem with a box domain $0.3\times 0.2$. %
The $N$-region, $\Omega_{N}$ with doping density $n_{d}$ consists of $x<0.1$ and $x>0.2$ and in the middle is the $P$-region $\Omega_{P}$ with density $n_a$, which is called N-P-N device.
$C(x)=n_{d}=10^{20}$ for $x<0.1$ or $x>0.2$ and $C(x)=-n_{a}=-6\times 10^{17}$ for $0.1<x<0.2$. %
Dirichlet boundary conditions $\varphi_{D}$, $n_{D}$, and $p_{D}$ are given on $x=0$ and $0.3$ and Neumann boundary conditions are given on other sides.
\par
\begin{figure}[tb]
  \centering
      \includegraphics[width=50mm]{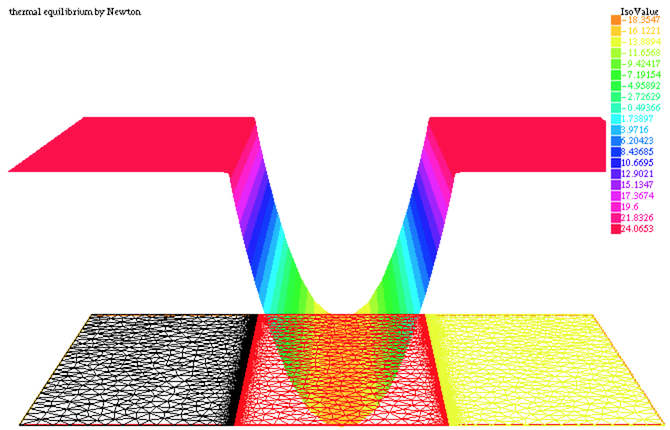}\ %    
\includegraphics[width=50mm]{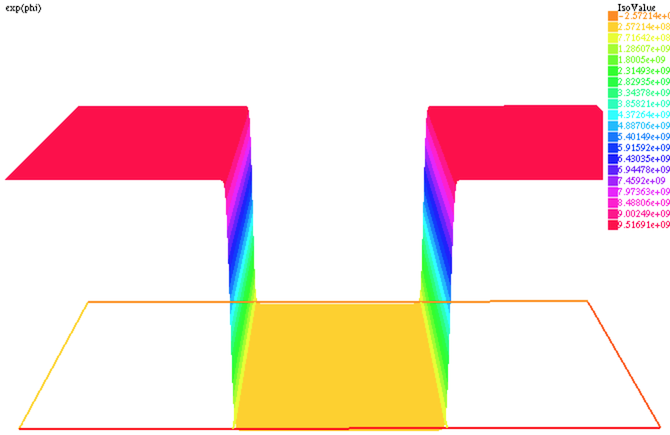}\ %
\includegraphics[width=50mm]{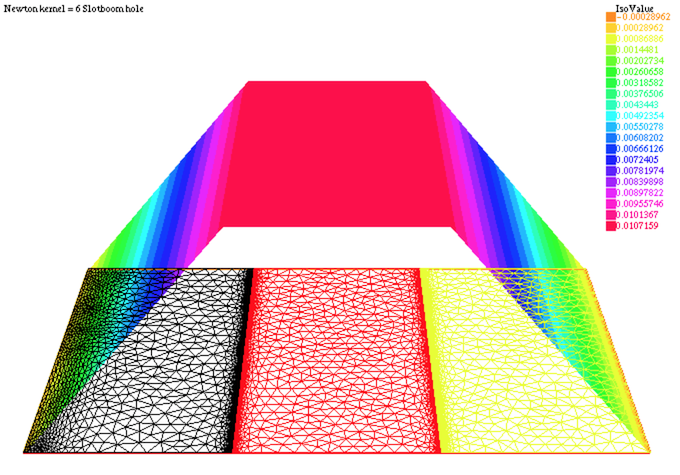}\\
  \caption{distribution of functions for N-P-N semiconductor device, left : electrostatic potential $\varphi_{\ast}$, middle : $e^{\varphi_{\ast}}$ in matrix coefficient, right : pseudo kernel of the hole $\xi^{0}$}
  \label{fig:npn}
\end{figure}

The electrostatic potential $\varphi_{\ast}$ by setting $\xi\equiv 1$ and $\eta\equiv 1$ and then satisfying $-\nabla(\lambda^{2}\nabla\varphi_{\ast})=e^{-\varphi_{\ast}}-e^{\varphi_{\ast}}+C(x)$ and $\varphi_{\ast}=\text{sinh}^{-1}(n_{d}/(2n_{i}))$ with $n_{i}=1.08\times 10^{10}$ on $x=0$ and $x=0.3$ is called thermal equilibrium. %
The left of Figure~\ref{fig:npn} shows distribution of $\varphi_{\ast}$ in the N-P-N device with $\Omega=(0, 0.3)\times(0, 0.2)$. %
A Newton iteration to obtain the thermal equilibrium is rather straightforward and is a part of the Gummel map~\cite{Brezzietal:2005}, which is a kind of fixed point method in total. %
To obtain a solution of the nonlinear system \eqref{eqn:drift-diffusion-nonlinear}, we will apply a Newton iteration starting from the thermal equilibrium $(\varphi, \eta, \xi)=(\varphi_{\ast}, 1, 1)$. %
By introducing expression on hole current density $J_{p}=-e^{-\varphi_{\ast}}\nabla\xi$, a mixed formulation of an elliptic equation with coefficient $e^{-\varphi_{\ast}}$ in the first step of the Newton step  is obtained as
\begin{equation}
  \int_{\Omega}e^{\varphi_{\ast}}J_{p}\cdot v-\int_{\Omega}\nabla\cdot v\,\xi
  - \int_{\Omega}\nabla\cdot J_{p}\,q=\int_{\Omega}f\cdot v-\int_{\Gamma_{D}}\xi_{D} v\cdot \nu,\label{eqn:mixed-Jp}
\end{equation}
where external force $f$ represents nonlinear coupling between electrostatic potential unknown $\varphi$ and hole unknowns $(J_{p}, \xi)$ and $\Gamma_{D}$ is a part of the Dirichlet boundary with $\xi_{D}=e^{\varphi_{D}}p_{D}$ and $\nu$ denotes the outer normal to boundary $\Gamma_{D}$. %
\par
Finite element matrix is obtained by discretization of a weak formulation of the weak formulation with RT0/P1 elements, where RT0 is the Raviart-Thomas finite element in the lowest order for vectorial unknown function $H(\text{div}\,;\,\Omega)=\{v\in L^{2}(\Omega)\,;\,\int_{\Omega}\nabla\cdot v <+\infty\}$ ~\cite{BoffiBrezziFortin:2013},
\begin{equation}
\overline{K} =
  \begin{bmatrix}
    A & B^{T} \\-B & 0
  \end{bmatrix}\label{eqn:stiffness-dd-hole}\,.
\end{equation}
Here mass and constraint matrices are defined using finite element basis functions $\{\varphi_{i}\}$ for $J_{p}$ and $\{\psi_{i}\}$ for $\xi$,
\begin{equation*}
  [A]_{i\,j}=\int_{\Omega} e^{\varphi_{\ast}} \varphi_{j}\cdot \varphi_{i},\quad
  [B]_{i\,j}=-\int_{\Omega}\nabla\cdot\varphi_{j}\psi_{i}\,.
\end{equation*}
We call this matrix as {\tt dd-hole}. %
The electrostatic potential $\varphi_{\ast}$ takes negative value in the $P$-region and ratio of $e^{\varphi_{\ast}}$ between $N$-region and $P$-region becomes below $10^{-16}$. By approximating $e^{\varphi_{\ast}}\simeq 0$ in $P$-region, the first term of \eqref{eqn:mixed-Jp} by the domain integration is replaced by
$\int_{\Omega_{N}}e^{\varphi_{\ast}}J_{p}\cdot v$  and
for arbitrary constant $c$, $(J_{p}^{0}, \xi^{0})$ that satisfies
\begin{equation*}
  J_{p}^{0}=0,\ \xi^{0}=c\text{ in }\Omega_{P}, \ -\nabla\cdot (e^{-\varphi_{\ast}}\nabla\xi^{0})=0 \text{ in }\Omega_{N}\text{ with }\xi^{0}=c \text{ on }\partial\Omega_{N}\cap\partial\Omega_{P}\text{ and }J_{p}^{0}=-e^{\varphi_{\ast}}\nabla\xi^{0}
\end{equation*}
will be the solution of the modified weak formulation. %
Shifted solution $(J_{p}, \xi)$ of \eqref{eqn:mixed-Jp} by $(J_{p}^{0}, \xi^{0})$ still almost satisfies the same weak formulation with difference as $\int_{\Omega_{P}}e^{\varphi_{\ast}}J_{p}\cdot v$, which is the residual of the approximation of $e^{\varphi_{\ast}}$ by zero.
This property confirms the stiffness matrix ${\widetilde K}$ is singular with one dimensional kernel, when all coefficients are stored in double precision. %
The middle and the right of Figure~\ref{fig:npn} show the exponential weight $e^{\varphi_{\ast}}$ with the thermal equilibrium and the kernel function $\xi^{0}$.
The maximum and minimum eigenvalues calculated by quadruple precision for given matrix in double precision, %
$\lambda_{\max}(K)=6.13488\times 10^{10}$, $\lambda_{\min}(K|_{\text{Im}(A)})=3.5704\times 10^{-12}$. %
The condition number of $K$ on $\text{Im}(A)$ is $\kappa(K|_{\text{Im}(K)})=1.7193\times 10^{22}$.
Applying double precision arithmetic with pivot threshold $\tau=0.01$, 13 entries are postponed and the Schur complement $S_{22}$, whose size is $13\times 13$, has one dimensional kernel. %
By this decomposition, $\lambda_{\max}(K_{11})=2.0080\times 10^{2}$, $\lambda_{\min}(K_{11})=3.6200\times 10^{-12}$,  and $\kappa(K_{11})=5.5472\times 10^{13}$.
\subsection{performance comparison of mixed  and pure precision arithmetic}
Here we summarize performance of the proposed algorithm in accuracy and computational speed. %
We used one core of Apple M1 Max CPU running at 3.23\,GHz, which has capability to perform two single precision arithmetic instead of one double precision in the same cycle.
For matrices {\tt phb1HYS} in Section~\ref{sec:matrix-market} and {\tt stokes} in Section~\ref{sec:numerical-Stokes}, double precision arithmetic is used for higher accuracy and single and double arithmetic are used as mixed precision which is labeled as mixed(double+single) in the Table~\ref{tab:performance}\,.
For matrix {\tt dd-hole} in Section~\ref{sec:semi-conductor}, quadruple precision arithmetic that is realized as double-double in QD library \cite{BaileyLiHida:2003} is used for higher accuracy and double and double-double arithmetic are used as mixed precision which is labeled as mixed(quadruple+double).
The error and the residual of the linear system is calculated from the RHS that is set to satisfy the solution is $[x]_{i}\equiv i (\text{mod}\ 11)$.
Since the later two matrices are singular, detected kernel dimension is also shown.
\begin{table}[h]
  \centering
  \caption{error, residual and elapsed time for factorization by pure and mixed precision arithmetic} 
  \begin{tabular}{cccc}
    \hline
    {\tt phb1HYS} & \multicolumn{3}{l}{$n=36,414$, $nnz=4,344,765$}\\
    &double &mixed(double+single) & single\\
    \hline
    error & $1.2650\times 10^{-6}$& $1.6647\times 10^{-6}$ & $3.1178\times 10^{-2}$ \\
    residual & $7.6201\times 10^{-16}$&$9.6442\times 10^{-16}$ & $4.4906\times 10^{-7}$\\
    time in second & 0.5328 & 0.4604 & 0.4053\\
    \hline
    \hline
    {\tt stokes} & \multicolumn{3}{l}{$n=374,536$, $nnz=21,146,848$}\\
    &double & mixed(double+single) & single\\
    \hline
    error & $4.3646\times 10^{-13}$& $2.0301\times 10^{-12}$ & $1.7046\times 10^{-2}$\\
    residual & $1.2730\times 10^{-15}$&$5.1881\times 10^{-15}$ & $7.0494\times 10^{-7}$\\
    dim. of kernel & 6 & 6 & 0\\
    time in second & 33.390 & 22.983 & 15.228\\
    \hline
    \hline
    {\tt hole} & \multicolumn{3}{l}{$n=40,323$, $nnz=401,243$}\\
    &quadruple &mixed(quadruple+double) & double\\
    \hline
    error & $5.3652\times 10^{-20}$& $1.3061\times 10^{-21}$ & $8.9382 \times 10^{-6}$\\
    residual & $5.3850\times 10^{-32}$&$1.5514\times 10^{-32}$ & $4.7860\times 10^{-16}$\\
    dim. of kernel & 1 & 1 & 1\\
    time in second & 16.599 & 2.6459 & 0.4064\\
    \hline
  \end{tabular}\label{tab:performance}
\end{table}
\section{Conclusions}\label{sec:5}
We have constructed a new hybrid algorithm for $LDU$-factorization for large sparse matrix introducing iterative solver for generation of Schur complement matrix in higher precision, where the matrix is decomposed into a union of moderate and hard parts. %
Numerical tests confirm the solution by the proposed algorithm by mixed precision arithmetic can keep accuracy as higher precision arithmetic. %
\par
When quadruple precision arithmetic are realized by using double-double data structure and are performed on the hardware equipped with fused multiply-add unit, ratio of arithmetic complexity of double-double to double is 25 to 1.
Therefore for the linear system that has huge condition number more than the range of the maximum floating digits of the double precision, mixed precision arithmetic with quadruple and double attains substantial speed-up, which was verified by a matrix from the semi-conductor problem.
\par
Since recent CPU has ratio of arithmetic complexity of double to float is 2 to 1, some speed-up is obtained, but solution phase by iteration procedure for recovering double precision accuracy for the Schur complement matrix masks the efficiency.
It is necessary to implement our hybrid factorization algorithm on the system with more single floating point arithmetic units than double and to evaluate the performance.
\par
For the solution phase of the large system with multiple RHSs to generate Schur complement, in the preconditioned part by forward/backward substitution in lower precision already well utilizes the BLAS level 3 routine, e.g., {\tt TRSM}, but
it is necessary to optimize SpMM operation in double and quadruple precision, because such kind of sparse linear algebra library is not provided yet.


\begin{thebibliography}{99}
  \bibitem{HighamMary:2022} Higham, N.~J. and Mary, T., Mixed precision algorithms in numerical linear algebra. \textit{Acta Numerica} (2022) 347--414.
  \bibitem{Geroge:1977} George, A. Numerical experiments using dissection methods to solve n by n grid problems. \textit{SIAM J. Numer. Analy.} (1977) \textbf{14}:161-–179. %DOI: 10.1137/0714011.
  \bibitem{SuzukiRoux:2014} Suzuki, A. and Roux, F.-X., A dissection solver with kernel detection for symmetric finite element matrices on shared memory computers. \textit{Int. J. Numer. Meth. Engng.} (2014) \textbf{100}:136--164.
  \bibitem{Saad:2003} Saad, Y., \textit{Iterative methods for sparse linear systems (2nd ed.)}. SIAM, (2003)
    \bibitem{WilliamsOlikerVuducShalfYelickDemmel} Williams, S., Oliker, L., Vuduc, R., Shalf, J., Yelick, K., and Demmel, J., Optimization of sparse matrix-vector multiplication on emerging multicore platforms, \textit{Parallel Computing} (2009) \textbf{35}:178--194.
  \bibitem{FarhatRoux:1994} Farhat, C, Roux F.-X., Implicit parallel processing in structural mechanics. \textit{Computational Mechanics Advances} (1994) \textbf{2}:1--124.
    \bibitem{Hecht:2002} Hecht, F., C++ tools to construct our user-level language., \textit{ESIAM: M2AN}, (2002) \textbf{36}:809--836.
    \bibitem{Brezzietal:2005} Brezzi, F., Marini, D.~D., Micheteletti, S., Pietra, P., Sacco, R. and Wang, S., Discretization of Semiconductor Device Problems(I), \textit{Handbook of numerical analysis}, Ciarlet, P.~G. ed, (2005) \textbf{13}:317--441.
    \bibitem{BoffiBrezziFortin:2013} Boffi, D., Brezzi, F. and Fortin, M. \textit{Mixed finite element methods and applications}, Springer, 2010.
      \bibitem{BaileyLiHida:2003} Bailey, D.~H., Li.~X.S., Hida, Y. QD: A double-double/quad-double package, Computer software, doi:10.11578/dc.20210416.14\ , (2003). % Jun.
\end{thebibliography}
\end{document}